\documentclass[12pt]{amsart}

\usepackage[T2A]{fontenc}
\usepackage[utf8]{inputenc}
\usepackage[english]{babel}

\usepackage{
 amsfonts
,amssymb
,amsmath
,mathabx
,colonequals
,amsthm
,yfonts
,fge
}%
\usepackage{a4}

\newcommand*{\op}{{\mathcal L}}

\newcommand*{\opB}{\op_B}

\newcommand*{\ppp}{{\mathfrak{p}}}
\newcommand*{\qqq}{{\mathfrak{q}}}
\newcommand*{\rrr}{{\mathfrak{r}}}
\newcommand*{\sss}{{\mathfrak{s}}}

\newcommand*{\pp}{{{p}}}
\newcommand*{\qq}{{{q}}}
\newcommand*{\rr}{{{r}}}
\newcommand*{\sS}{{{s}}}

\newcommand*{\argz}{{}}

\newcommand*{\llll}{{\ell}}

\newcommand*{\delTa}{{\mathfrak D}\!}

\newcommand{\half}{\mbox{\scriptsize$1\over2$}}
\newcommand{\mhalf}{\mbox{\scriptsize{-$1\over2$}}}
\newcommand{\halfi}{\mbox{\scriptsize$\Imi\over2$}}
\newcommand{\mhalfi}{\mbox{\scriptsize{-$\Imi\over2$}}}

\newcommand{\quarti}{\mbox{\scriptsize$\Imi\over4$}}
\newcommand{\mquarti}{\mbox{\scriptsize$-\Imi\over4$}}

\newtheorem{theorem}{Theorem}
\newtheorem{proposition}[theorem]{Proposition}

\newtheorem{definition}{Definition}

\newcommand*{\Eqs}{{Eq.s}}

\newcommand{\Imi}{\mathrm{i}}

\newcommand{\eiphi}{{\Phi}}

\newcommand{\calF}{{\Theta}}
\newcommand{\CCcalF}{{\tilde\Theta}}

\newcommand{\eP}{{\Psi}}

\newcommand{\somEE}{    E    }
\newcommand{\EEpm}[1]{E_{\{\! #1 \!\}}}

\newcommand{\backSlash}{\backslash}
\renewcommand*{\backSlash}{\,\fgebackslash\,}

\newcommand{\covered}[1]{\mbox{\textswab{#1}}}
\newcommand{\domain}{\covered{C}}
\newcommand{\subdomain}{\mbox{$\mathstrut^\backprime{}\covered{C}$}}

\newcommand{\diag}{\,{\mathrm{diag}}}

\newcommand{\mono}{\mathcal{M}}
\newcommand{\monodrom}{monodromy}

\newcommand{\AAACCC}{u}
\newcommand{\BBBDDD}{v}
\newcommand{\KKKLLL}{w}

\newcommand{\AAA}{\AAACCC_{\!+}}
\newcommand{\CCC}{\AAACCC_{\!-}}
\newcommand{\BBB}{\BBBDDD_{\!+}}
\newcommand{\DDD}{\BBBDDD_{\!-}}
\newcommand{\KKK}{\KKKLLL_{\!+}}
\newcommand{\LLL}{\KKKLLL_{\!-}}

\newcommand{\puncturedS}{^\backprime{\!}S^1}

\begin{document}

\title[Monodromy map]%
{
Square root of the monodromy map
for the equation of RSJ model of Josephson junction
}

\author{Sergey I. Tertychniy
}%

\address%
 {
 Russian Metrological Institute of Technical Physics and Radio Engineering
 (VNIIFTRI), Mendeleevo, 141570, Russia}%
 \thanks{
 Supported in part by RFBR grant N 17-01-00192.
 }
\begin{abstract}\noindent
Several noteworthy  properties of  the differential equation utilized
for the modeling of Josephson junctions are considered. The explicit
representation of the monodromy transform of the space of its
solutions is given. In case of positive integer order, the
transformation interpreted as the square root of the monodromy
transformation noted is derived making use of 
a symmetry the associated linear second order differential equation possesses.

\end{abstract}

\maketitle

The present notes
are %
devoted to discussion of
some
noteworthy properties of 
the differential equation
\begin{equation}
   \label{eq:010}
\dot\varphi +\sin\varphi=B+A\cos\omega t,
\end{equation}
in which  $\varphi=\varphi(t) $ is the unknown function,
the symbols
$A,B,\omega$ stand for some real constants, and $t$ is a free real
variable,
the %
dot %
denoting %
derivation
with respect to $t$.
Eq.{~}\eqref{eq:010} and its generalizations
are of interest, in particular, in view of their application
in
a number of models in physics,
mechanics, and geometry \cite{Foo,FLT}.
Most widely Eq.{~}\eqref{eq:010} is known as the equation
utilized in the so called
RSJ model of Josephson junction
\cite{
St,McC,Ba,MK,Sch}
which applies
if 
the effect
of the junction 
electric capacitance 
is negligible (the case of so called overdamped Josephson junctions).

Eq.{~}\eqref{eq:010} is equivalent to the 
Riccati equation
\begin{eqnarray}
          \label{eq:020}
&&
 {\eiphi}'
=(2\,\Imi\,\omega\,z)^{-1}(1-{\eiphi}^2)+
(\llll\, z^{-1}+\mu(1+z^{-2})){\eiphi},
\end{eqnarray}
where $ {\eiphi}= {\eiphi}(z)$
is a  
holomorphic function of the free complex variable $z$,
the prime denoting the derivative with respect to the latter, 
and $\llll,\mu,\omega $ are the constant parameters.
Indeed, the formal substitutions
\begin{equation}
\label{eq:030}
z \leftleftharpoons 
e^{\Imi\omega t},\; {\eiphi}(z)
\leftleftharpoons 
e^{\Imi\varphi(t)}
\end{equation}
convert Eq.{~}\eqref{eq:020}
to Eq.{~}\eqref{eq:010} get with the parameters
related to
the
parameters 
involved in the former equation
by the transformation 
\begin{equation}
 \label{eq:040}
A=2\omega\mu, \; B=\omega\llll.
\end{equation}

The nonlinear equation \eqref{eq:020} has the only singular point $z=0$
but its solutions may diverge at
 of some other values of argument.
Nevertheless, one can claim the following \cite{A6}:
\begin{proposition}
Let
the constants $\llll,\mu$, and $\omega>0$ be real and let
$ {\eiphi}(z)$ be a solution to Eq.{~}\eqref{eq:020}
holomorphic at $z=1$
such that
$ |{\eiphi}(1)|=1  $. Then $ {\eiphi}(z)$  is also holomorphic
in some vicinity of the curve $|z|=1$; moreover, if  $|z|=1$
then $ |{\eiphi}(z)|=1  $.
\end{proposition}
Indeed, any solution to Eq.{~}\eqref{eq:010} is real analytic and
can be extended to  the whole real axis $\mathbb{R}$.
Carrying out
analytic continuation of $\varphi(t)$ from $\mathbb{R}$
to some its vicinity in  $\mathbb{C}$, one obtains the function
${\eiphi}=\exp(\Imi\varphi(t))$ which is holomorphic in $t$
and which
can also be considered as a holomorphic function of
the variable
$z=\exp(\Imi\omega t)$
varying in some open set embodying
the curve %
 $|z|=1$. The function $ {\eiphi}(z)$ possesses
the properties asserted above, obviously.

It is also obvious that if we are given the function
$ {\eiphi}(z)$ obeying Eq.{~}\eqref{eq:020} and unimodular on
the curve $|z|=1$, then the real-valued smooth function $\varphi(t)$
such that $ e^{\Imi \varphi(t)}= {\eiphi}(e^{\Imi\omega t})  $
can be constructed.
It verifies
Eq.{~}\eqref{eq:010}, evidently.

It is in order now to comment on
the term ``the curve $|z|=1$''
used above
instead of
something like
 ``the unit circle $S^1\subset\mathbf{C}$''
which one could argue to be more customary. The point is that,
strictly speaking,
apart of the very special  conditions, 
a solution to Eq.{~}\eqref{eq:020} can not be holomorphic on $S^1$.
The rationale is here fairly
simple: indeed, there is no reason why
a generic solution $\varphi(t)$ to
Eq.{~}\eqref{eq:010} should  obey the constraint
$\varphi(\half T)=\varphi(\mhalf T)\, (\!\!\!\mod2\pi)$, where
$$T=2\pi\omega ^{-1}$$
is the period of the right hand side expression in Eq.{~}\eqref{eq:010}.
Accordingly, following the way of constructing of the function $ {\eiphi}$
via the function $\varphi$ utilized above, one obtains 
$ {\eiphi}(-1+\!0)-{\eiphi}(-1-\!0)=e^{\Imi\varphi(\half{}T)}-e^{\Imi\varphi(\mhalf{}T)}\not=0$,
meaning that  $ \eiphi $, when considered on $S^1$, proves to be not continuous at
$-1$.

 There is, definitely, no singularity of $ {\eiphi}$ at $-1$ and the
 deficiency 
in the above construction originates in the  improper selection of the $ {\eiphi}$ domain.
Generally speaking, it can not be the complex plane or any its subset; instead,
the universal cover \domain{}  of the punctured complex plane
$\mathbb{C}^*=\mathbb{C}\backSlash 0$ (with the 
subset of the isolated
singular
points of $ {\eiphi}$ removed)
has to be utilized.
In this setting,
 the image of
the ``$t$-axis'', produced  by the map extending
the transformation \eqref{eq:030} to \domain, is not  $S^1$
but the non-compact curve covering $S^1$.
It is this curve which, admitting some abuse of notations,
 was referred to  as ``the curve $|z|=1$''.

Here we shall not, however, consider solutions to Eq.{~}\eqref{eq:020} on their
whole domains  but only on the sub-domain \subdomain{} 
which is in bijective correspondence with (projects to)
the subset ${^\backprime{}\mathbb{C}^*}$ of
$\mathbb{C}^*$ obtained by removal of the ray of negative real numbers,
 ${^\backprime{}\mathbb{C}^*}=\mathbb{C}\backSlash\mathbb{R}_{\le0}$, (and
removal, for each ${\eiphi}$, its singular points, if any).
In most cases, \subdomain{}  can be (and will be) considered
to be undistinguished from ${^\backprime{}\mathbb{C}^*}$.
However,
 there are two boundaries of \subdomain,
which projects to the same (removed) ray $\mathbb{R}_{<0}$, on which the values of
$ \eiphi $ do not coincide.
To mirror such a difference,
we may identify these boundaries with the two edges
of the corresponding cut in ${\mathbb{C}^*}$.
It is then convenient to refer to the points
of the cut edge contacting the half-plane
$\{z\in \mathbb{C}^*,\Im z>0\}  $ by the symbol $\rho\, e^{\Imi\pi}  $, where $\rho$
stands for a positive real number, $\rho\in\mathbb{R}_{>0}$, and
by the symbol $\rho\, e^{-\Imi\pi} $ for a point belonging to the cut edge contacting the
half-plane
$\{z\in \mathbb{C}^*,\Im z<0\}  $. If $\rho=1$ the factor $\rho$  is omitted.
In such a framework, 
 the portion of $S^1$ coming to be in ${^\backprime{}\mathbb{C}^*}$
is ``the punctured circle''
\begin{equation}
\label{eq:050}
\puncturedS=\{z\in\mathbb{C}, |z|=1, z\not=-1\}.
\end{equation}%
It approaches near the ends the (distinct)
boundary points {\em denoted\/}  $e^{\Imi\pi}$ and $e^{-\Imi\pi}$.
Each solution  to Eq.{~}\eqref{eq:020} is holomorphic and non-zero
in some vicinity of $ \puncturedS $.
Besides, at the boundary points of $ \puncturedS $,
it holds
\begin{equation}
\label{eq:060}
{\eiphi}( e^{\Imi\pi} )=e^{\Imi\varphi(\half{}T)}
 \;(\not=)\;
{\eiphi}( e^{-\Imi\pi} )
=e^{\Imi\varphi(\mhalf{}T)}.
\end{equation}

The following statement is an obvious consequence of the
periodicity of the right hand side of Eq.{~}\eqref{eq:010}:
        \begin{proposition}
Let the function $\varphi(t) $ verify Eq.{~}\eqref{eq:010}.
Then the function $\varphi_M(t)\colonequals \varphi(t+T)$ is also
a solution to Eq.{~}\eqref{eq:010}.
     \end{proposition}
Let the solution
 ${\eiphi}(z)$ to Eq.{~}\eqref{eq:020}
 be constructed in accordance with the algorithm
specified above. The above statement and definitions imply the following.
\begin{proposition}
There exists the solution ${\eiphi}_M(z)$ to Eq.{~}\eqref{eq:020}
such that
\begin{equation}
\label{eq:070}
                  {\eiphi}_M( e^{-\Imi\pi})={\eiphi}(e^{\Imi\pi}).
\end{equation}
Moreover, if \eqref{eq:070} holds true
for some solutions ${\eiphi}$, ${\eiphi}_M$ 
 then
\begin{equation}
{\eiphi}_M(\rho\, e^{-\Imi\pi})={\eiphi}(\rho\,e^{\Imi\pi}) 
\end{equation}
for all $\rho>0$ excluding ones  
at which ${\eiphi}(\rho\,e^{\Imi\pi})$ is not
analytic (i.e.\ is undefined).
\end{proposition}
Indeed, the function $ {\eiphi}_M(z) $ wanted can be constructed
from the function $\varphi_M(t)$ in the same way as the function
$ {\eiphi}(z) $ is constructed from
$\varphi(t)$. The function $\varphi_M(t)$ is actually
some ``portion'' of the maximally extended solution $\varphi(t)$
get from  the segment $(\mhalf{}T+T,\half{}T+T)$
and ``put down'' 
to the segment $(\mhalf{}T,\half{}T)$
considered as the common domain with $\varphi(t)$;
 similarly,
the function  ${\eiphi}_M$ is, essentially,
the function ${\eiphi}$ get on the sub-domain adjacent via the common boundary
with
\subdomain{}  and considered 
on the
same domain
with ${\eiphi}$
 trough their projections to the common area
 ${^\backprime{}\mathbb{C}^*}$
which can be considered equivalent to \subdomain.

Alternatively,
the function  ${\eiphi}_M$ can also be defined as the result of
 point-wise analytic continuations in ${\mathbb{C}^*}$
of the function ${\eiphi}$
defined on ${^\backprime{}\mathbb{C}^*}$
along the
full circles with centers situated at zero
which are
passed in
the
counter-clockwise direction (or along the curves avoiding
${\eiphi}$ singularities and homotopic to such circles).
The latter interpretation allows one to refer to 
transformation
$M\!\!:{\eiphi}\mapsto{\eiphi}_M  $ as the {\em \monodrom{}  map\/}
which acts  
on the space 
of solutions to Eq.{~}\eqref{eq:020}.

We are now ready to formulate the first non-evident 
result the present notes are devoted to.
\begin{theorem}\it
Let a solution $\eiphi $
to Eq.{~}\eqref{eq:020}
holomorphic
in
some vicinity of
$\puncturedS$ be given.
Let also $\eP=\eP(z) $ be a  solution of the linear homogeneous
first order ordinary differential equation
\begin{equation}
                 \label{eq:080}
2\Imi\omega z\eP'=( {\eiphi} + {\eiphi}^{-1})\eP. 
\end{equation}
Let, finally,
\begin{equation}
                  \label{eq:090}
\eP(1)=1 \mbox{ and }
| \eiphi(1) |=1.  
\end{equation}
Then the formula
\begin{equation}
              \label{eq:100}
\begin{aligned}
\eiphi_{M}(z)=&
\left(
e^{\half{}P(\half T)}
\cos\half\varphi(\half{}T)
\cdot\eP(z)^{\half}\eiphi(z)^{\half} \right.
\\&\left.
+\Imi
e^{\half{}P(\mhalf T)}
\sin\half(\varphi(\half{}T)-\varphi(\mhalf{}T))
\cdot\eP(1/z)^{\half}\eiphi(1/z)^{\mhalf}
\right)\times
\\&
\left(
e^{\half{}P(\half T)}
\cos\half\varphi(\half{}T)
\cdot\eP(z)^{\half}\eiphi(z)^{\mhalf} \right.
\\&\left.
-\Imi
e^{\half{}P(\mhalf T)}
\sin\half(\varphi(\half{}T)-\varphi(\mhalf{}T))
\cdot\eP(1/z)^{\half}\eiphi(1/z)^{\half}
\right)^{-1}\hspace{-1em}
\end{aligned}
\end{equation}
in which the  
continuous (and then necessarily real analytic) function $ \varphi $
is determined by the equation ${\eiphi}(e^{\Imi\omega{}t})=e^{\Imi\varphi(t)}, t\in(-\half{}T,\half{}T)  $,
yields the explicit representation of the result of the
 \monodrom{}  transformation of the function ${\eiphi}$.
\end{theorem}
The above assertion means that any solution to Eq.{~}\eqref{eq:020}
can be extended from its sub-domain 
with the closure equal to 
\subdomain{}
to the whole domain with the closure equal to
\domain{}
by means of certain algebraic
transformations 
(provided the function
 $\eP$ had once only been computed 
 on
\subdomain{}).

Before proving these, it is
worth commenting on the existence of 
$\eP$. Given $ {\eiphi} $,
it is determined
on the base of the equality
$$
{\eP}(e^{\Imi\omega{}t})=e^{P(t)}, t\in(-\half{}T,\half{}T),\;\mbox{where}\;
P(t)=\mbox{$\int^t_0\cos\varphi(\tilde t)\,d\,\tilde t$},
$$
reducing, therefore, to the quadrature and subsequent analytic
continuation of its result
from an arc of the curve $|z|=1$. Notice also that for such $\eP $
the square root
$\eP^{1/2}$ is uniquelly defined via analytic continuation of $
e^{P(t)/2} $. The functions $\eiphi^{\pm1/2}$ are endowed with  unique
values in a similar way.

The proof of the formula \eqref{eq:100} splits 
into  two steps.
First, its right hand side is evaluated for the argument $z=e^{-\Imi\pi}$.
Performing substitutions in accord with definitions,
 one obtains
$ {\eiphi}_M(e^{-\Imi\pi})=e^{\varphi(\half{}T)}={\eiphi}( e^{\Imi\pi} )  $.
Second, the expression \eqref{eq:100} is substituted
into Eq.{~}\eqref{eq:020}. Then, upon elimination of the derivatives
$ {\eiphi}' $ and $ {\eP}' $ with the help of Eq.{~}\eqref{eq:020}
and Eq.{~}\eqref{eq:080}, respectively, the identical equality follows.
Thus, the expression \eqref{eq:100} verifies the
first order differential equation
\eqref{eq:020}
and obeys the initial condition \eqref{eq:070} which
distinguishes the solution representing the
\monodrom{}  transformation of $ {\eiphi} $.
The identical coincidence follows and we are done.

It will be further assumed throughout that the parameter $\llll$ is a positive integer,
$$
\llll\in\mathbb{N}.
$$
We set up 
the following definition \cite{BT1}.
\begin{definition}
Let the four sequences $p_k, q_k, r_k, s_k,\; k=0,1,2\dots $ of functions
of the complex variable $z$ and the constant parameters $\llll,\mu,\lambda=(2\omega)^{-2}-\mu^2  $
be
defined by means of the following recurrent scheme
\begin{eqnarray}
&&\begin{aligned}\hspace{0.9em}
\pp_{0}=0,\; \qq_{0}=1,\;\rr_{0}=z^{-2},\; \sS_{0}=-\mu;
\end{aligned}
\\
&&
\left\{
\begin{aligned}
\pp_{k}=&\:
                (1 - \llll)z\,\pp_{k-1}  + \qq_{k-1} + z^2 \pp_{k-1}', \\
\qq_{k}=&\:
               z^2 (-\lambda  + (\llll+1) \mu z) \pp_{k-1} +  \mu\left(1- z^2\right) \qq_{k-1}
               +  z^2  \qq_{k-1}';
\\[0.2em]
    \rr_{k}=&\:
2 (k-2)z\,\rr_{k-1} - \sS_{k-1} -  z^2  \rr_{k-1}',
\\
     \sS_{k}=&\:
z^2 \left(\lambda-\left(\llll + 1\right) \mu z\right)\rr_{k-1}
+ \left(\left(2k-\llll-3\right) z + \mu\left(z^2-1\right) \right)\sS_{k-1}
- z^2  \sS_{k-1}',
\end{aligned}
\right.
\end{eqnarray}
We pick up
their 
``diagonal''
representatives denoting them
$\ppp,\qqq,\rrr,\sss$, i.e.\ define
\begin{equation}
\ppp= \pp_{\llll}, \qqq=\qq_{\llll}, \rrr=\rr_{\llll}, \sss=\sS_{\llll}.
\end{equation}
\end{definition}
It can be shown that the functions $\ppp,\qqq,\rrr,\sss$
are the polynomials in $z$ of the degrees $2\llll-2, 2\llll,
2\llll-2,2\llll$,
respectively \cite{BT1};
 they are polynomial in the parameters $\lambda$  and $\mu$ as well.

We define now the following four holomorphic functions
 $\eiphi_B, \eP_B, \calF_B , \CCcalF_B $
representing them 
in terms of the two other holomorphic functions
$\eiphi$ and $\eP$ of a complex variable
and the functions $\varphi, P $ of a real variable,
the latter pair being  evaluated
for several fixed values of their argument alone:
\begin{equation}\label{eq:130}
\begin{aligned}
\eiphi_B(z)=&
-
\left(
2\Imi e^{\half P(\half T)}
(\delTa_{+ }  \LLL \CCC + \delTa_{- }  \KKK \AAA  )
\cdot\eP(z)^{\half}\eiphi(z)^{\half}
\right.\\&\hphantom{-}
+
\big(
-\delTa_{+ }  \LLL
 (e^{\half P(\half T)} \CCC + e^{\half P(\mhalf T)} \DDD)
\\&\hphantom{-}\left.\hspace{1.7em}
+\, \delTa_{- }  \KKK
 (e^{\half P(\half T)}  \AAA  + e^{\half P(\mhalf T)} \BBB)
\big)
\cdot\eP(1/z)^{\half}\eiphi(1/z)^{\mhalf}
\right)\times
\\&
\left(
-2\Imi e^{\half P(\half T)}
(\delTa_{+ }  \LLL \CCC + \delTa_{- }  \KKK \AAA  )
\cdot\eP(z)^{\half}\eiphi(z)^{\mhalf}
\right.\\&
+
\big(
-\delTa_{+}  \LLL
 (e^{\half P(\half T)} \CCC + e^{\half P(\mhalf T)} \DDD)
\\&\left.\hspace{1.7em}
+\delTa_{- }  \KKK
 (e^{\half P(\half T)}  \AAA  + e^{\half P(\mhalf T)} \BBB)
\big)
\cdot\eP(1/z)^{\half}\eiphi(1/z)^{\half}
\right)^{-1},
\end{aligned}
\end{equation}
\begin{equation}
                    \label{eq:140}
\eP_B(z)=(2\Imi)^{-1}(\calF_B(z)-\CCcalF_B(z)), \mbox{ where}
\end{equation}
\begin{equation}
                   \label{eq:150}
\begin{aligned}
\calF_B(z)=&
       (\cos\varphi(0))^{-1}\times
\\&\hspace{-5ex}
\left(
-\Imi\big(
\delTa_{-}\KKK
(
(2\sin\varphi(0)-1)\,
 e^{\half P(\half T)} \AAA
-e^{\half P(\mhalf T)} \BBB
)
\right.\\&\hspace{-5ex}\hspace{1.0em}
+
\delTa_{+}\LLL
(
(2\sin\varphi(0)+1)\,
 e^{\half P(\half T)} \CCC
+e^{\half P(\mhalf T)} \DDD
)
\big)
          \cdot\eP(z)^{\half}\eiphi(z)^{\mhalf}
\\&\hspace{-5ex}
+\big(
     - \delTa_{-}\KKK
            (
              (\sin\varphi(0)-2)\,
                e^{\half P(\half T)} \AAA
                 +\sin\varphi(0)\,
                  e^{\half P(\mhalf T)} \BBB
             )
\\&\left.\hspace{-0.3em}
+
 \delTa_{+}\LLL
(
(\sin\varphi(0)+2)\,
 e^{\half P(\half T)} \CCC
+\sin\varphi(0)\,
e^{\half P(\mhalf T)} \DDD
)
\big)\times 
\right.\\&\left.\hspace{18em}
          \eP(1/z)^{\half}\eiphi(1/z)^{\half}
\right)\times
\\&
\left(
-2\Imi e^{\half P(\half T)}
(\delTa_{+ }  \LLL \CCC + \delTa_{- }  \KKK \AAA  )
\cdot\eP(z)^{\half}\eiphi(z)^{\mhalf}
\right.\\&
+
\big(
-\delTa_{+}  \LLL
 (e^{\half P(\half T)} \CCC + e^{\half P(\mhalf T)} \DDD)
\\&\left.\hspace{1.9em}
+\delTa_{- }  \KKK
 (e^{\half P(\half T)}  \AAA  + e^{\half P(\mhalf T)} \BBB)
\big)
\cdot\eP(1/z)^{\half}\eiphi(1/z)^{\half}
\right)^{-1},
\end{aligned}
\end{equation}
\begin{equation}
                               \label{eq:160}
\begin{aligned}
\tilde\calF_B(z)=&
       (\cos\varphi(0))^{-1}\times
\\&\hspace{-5ex}
\left(
\Imi\big(
            -\delTa_{-}\KKK
(
-(2\sin\varphi(0)-1)\,
  e^{\half P(\half T)}
            \AAA
 +e^{\half P(\mhalf T)}
            \BBB
)
\right.\\&\hspace{-5ex}\hspace{1.0em}
+
\delTa_{+}\LLL
(
(2\sin\varphi(0)+1)\,
  e^{\half P(\half T)} \CCC
 +e^{\half P(\mhalf T)} \DDD
)
\big)
          \cdot\eP(z)^{\half}\eiphi(z)^{\half}
\\&\hspace{-5ex}
+\big(
     - \delTa_{-}\KKK
            (
              (\sin\varphi(0)-2)\,
                 e^{\half P(\half T)} \AAA
                 +\sin\varphi(0)\,
                  e^{\half P(\mhalf T)} \BBB
             )
\\&\left.\hspace{-0.3em}
+
 \delTa_{+}\LLL
              (
                (\sin\varphi(0)+2)\,
                  e^{\half P(\half T)} \CCC
                 +\sin\varphi(0)\,
                  e^{\half P(\mhalf T)} \DDD
)
\big)\times 
\right.\\&\left.\hspace{18em}
          \eP(1/z)^{\half}\eiphi(1/z)^{\mhalf}
\right)\times
\\&
\left(
2\Imi e^{\half P(\half T)}
    (\delTa_{+ }  \LLL \CCC
   + \delTa_{- }  \KKK \AAA  )
\cdot\eP(z)^{\half}\eiphi(z)^{\half}
\right.\\&
+
\big(
-\delTa_{+}  \LLL
 (e^{\half P(\half T)} \CCC + e^{\half P(\mhalf T)} \DDD)
\\&\left.\hspace{1.9em}
+
  \delTa_{- }  \KKK
  (e^{\half P(\half T)}  \AAA
 + e^{\half P(\mhalf T)} \BBB)
\big)
\cdot\eP(1/z)^{\half}\eiphi(1/z)^{\mhalf}
\right)^{-1}.
\end{aligned}
\end{equation}
Above,
the following coefficient shortcuts
\begin{eqnarray}
\AAACCC_{\pm} &\colonequals&
(-1)^{\llll}e^{\halfi\varphi(\half T)}
\pm\Imi\, e^{\mhalfi\varphi(\half T)}
\nonumber \\
\BBBDDD_{\pm} &\colonequals&
e^{\halfi\varphi(\mhalf T)}
\pm
\Imi\,(-1)^{\llll}  e^{\mhalfi\varphi(\mhalf T)},
\nonumber\\
\KKKLLL_{\pm}&\colonequals& e^{\halfi\varphi(0)} \pm \Imi\, e^{\mhalfi\varphi(0)},
\nonumber \\
                   \label{eq:170}
\delTa_{\pm}&\colonequals&\ppp(1)\pm2\omega\rrr(1),
\end{eqnarray}
are utilized.

It is worth noting that 
the involvement
of the
functions $\varphi, P $
in
 \Eqs{~}\eqref{eq:130}-\eqref{eq:160}
is not obligatory.
Their values utilized there 
can be expressed
in terms of the functions $\eiphi,\eP $ alone,
provided the following identifications
 are taken into account ({\it cf} \Eqs{~}\eqref{eq:060}):
\begin{eqnarray}
                   \label{eq:180}
&&e^{\Imi\varphi(\half{}T)}=
{\eiphi}( e^{\Imi\pi} ),
 \;
e^{\Imi\varphi(0)}=
{\eiphi}(1 ),
 \;
e^{\Imi\varphi(\mhalf{}T)}
=
{\eiphi}( e^{-\Imi\pi} );
\;
\\
                   \label{eq:190}
&&
e^{P(\half{}T)}
=
{\eP}( e^{\Imi\pi} ),
\;
e^{P(\mhalf{}T)}
=
{\eP}( e^{-\Imi\pi} ).
\end{eqnarray}


The (second) non-obvious result
to be here reported
is as follows: 
                   \begin{theorem}\it
Let the functions ${\eiphi}$ and  $\eP $
verify the equations \eqref{eq:020} and \eqref{eq:080}, respectively,
obeying also the constraints \eqref{eq:090}.
Then
\begin{itemize}
\item the functions 
$\eiphi_B$ \eqref{eq:130} and $\eP_B$ \eqref{eq:140}
verify  the same equations and constraints as ${\eiphi}$ and  $\eP $,
respectively;
\item
the transformation $B: (\eiphi,\eP) \mapsto (\eiphi_B,\eP_B) $ 
repeated twice coincides with the \monodrom{}  transformation $M$.
\end{itemize}
\end{theorem}
Thus, $B$ can be considered as a square root of $M$.

The first assertion is proven by straightforward computation.
With regard to
the second one, 
 we replace here its formal proof  with  outline
of derivation of the very formulas \eqref{eq:130}-\eqref{eq:160}
demonstrating how they had been arisen. Besides, along the way, a  
profound 
relationship of the equation \eqref{eq:020} (and \eqref{eq:010})
with another family of differential equations is
demonstrated.

To that end, let us consider the two holomorphic functions
$\EEpm{\pm}(z)$
defined through the functions $\eiphi, \eP  $ as follows:
\begin{equation}
                  \label{eq:200}
\begin{aligned}
E_{\{\pm\}}(z)&
\colonequals 
2^{-1}
e^{\mu(z+1/z-2)/2}
z^{-\llll/2}
\times\\&\hspace{2ex}
\left\lgroup
\frac{1\pm\Imi}{\sqrt2}
(\eP(z)\eiphi(z))^{1/2}
+
\frac{1\mp\Imi}{\sqrt2}
(\eP(1/z)/\eiphi(1/z))^{1/2}
\right\rgroup.
\end{aligned}
\end{equation}
Straightforward calculation proves
the following equalities
\begin{equation}
                        \label{eq:210}
\EEpm{\pm}'(z)=\pm(2\omega)^{-1} z^{-\llll-1}   \EEpm{\pm}(1/z)
+\mu \EEpm{\pm}(z),
\end{equation}
taking place
provided the functions 
$\eiphi, \eP  $
obey the equations 
\eqref{eq:020} and \eqref{eq:080}, respectively.
\Eqs{} \eqref{eq:210} imply, in turn,
the fulfillment of the equation 
\begin{equation}
           \label{eq:220}
z^2 \somEE''(z)
+\big((\llll+1)
 z
+ \mu (1-z^2) \big) \somEE'(z)
+\big
(  -\mu  (\llll+1) z
+
\lambda
 \big
)
\somEE(z)=0
\end{equation}
by the both functions $\somEE=\EEpm{+}$
and $\somEE=\EEpm{-}$.

The equations of the form
\eqref{eq:220} with arbitrary constant parameters
$\llll, \lambda, \mu $ constitute a subfamily
of the family of so called double confluent Heun equations,
see Refs.~\cite{SW,SL,HP}.

Eq.{~}\eqref{eq:220} is a linear homogeneous differential equation with
coefficients holomorphic everywhere except zero.
Hence their solutions, including $\EEpm{\pm}$, are  holomorphic
everywhere except zero including the points of divergence
and roots
of
the solution ${\eiphi}$ of the non-linear equation \eqref{eq:080}
connected with $\EEpm{\pm}$ via \Eqs{~}\eqref{eq:200}.
At the same time, the common singular point $z=0 $ for all the functions
$ \EEpm{+},\EEpm{-} $ and $ {\eiphi}, \eP $ is actually the branching point
of their common domain, the universal cover \domain{}  of
the punctured complex plane
$\mathbb{C}^*$ (for $ {\eiphi}$ and  $\eP $, with their singular
points removed).
Being defined on \domain,
the functions $ \EEpm{\pm}$
behave like multi-valued functions on $\mathbb{C}^*$
and
may thus undergone the
\monodrom{}  transformation without violation of fulfillment of Eq.~\eqref{eq:220}.
Similarly to the case of 
 solutions to Eq.~\eqref{eq:020},
the
\monodrom{}  transformation of $ \EEpm{\pm}$
can be understood as 
point-wise analytic continuations
along the arcs projected to
full circles with centers situated at zero
which are passed in
the
counter-clockwise direction
(as opposed to the case of $\eiphi$,
no singular points can now be encountered on such arcs).

Substituting 
 $z=1$ into  \eqref{eq:200}, one gets
\begin{equation}
            \label{eq:230}
\mbox{$
E_{\{\pm\}}(1)=\mp\sin(\half 
(\varphi(0)\mp\half\pi)).  
$}
\end{equation}
Thus, if $\varphi(0)\not=\half\pi\,(\hspace{-1.7ex}\mod \pi)$ (i.e.\ if
\begin{equation}
\eiphi(1)^2\not=-1,
\end{equation}
 see \Eqs{~}\eqref{eq:180}) then 
$\EEpm{+}(z)\not\equiv0\not\equiv\EEpm{-}(z)$ 
and the functions
$\EEpm{+}$ and $\EEpm{-}$ are linear independent.
Moreover, since the linear space of solutions to Eq.{~}\eqref{eq:220}
is two-dimensional, the functions $\EEpm{\pm}$ constitute
its basis and any solution to Eq.{~}\eqref{eq:220}
can be represented as their linear combination with constant coefficients.
Thus, the two formulas \eqref{eq:200} ensure,
in fact,
the explicit
representation of all the solutions to Eq.{~}\eqref{eq:220}
in terms of any generic solution ${\eiphi}  $ to Eq.{~}\eqref{eq:020}
and some related quadrature (the function $\eP$).

Conversely, the formula
\begin{eqnarray}
              \label{eq:240}
\hspace{2em}
{\eiphi}^{(\alpha)}(z)
\!\!   &\colonequals&\!\!
-\Imi  z^l
\frac{
\cos(\half 
\alpha)\EEpm{ +}(z)
+\Imi
\sin(\half 
\alpha)\EEpm{ -}(z)
}{
\cos(\half 
\alpha)\EEpm{ +}(1/z)
-\Imi
\sin(\half 
\alpha)\EEpm{ -}(1/z)
}
\end{eqnarray}
in which $\alpha$ stands for an arbitrary real number,
represent a solution to
Eq.{~}\eqref{eq:020} obeying the constraint $|\eiphi(e^{\Imi\omega{}t})|=1$,
provided the functions $\EEpm{\pm}$ obey
\Eqs{~}\eqref{eq:210} and $\Im\EEpm{\pm}(1)=0$.
The composition of the transformations
Eq.{~}\eqref{eq:200} and Eq.{~}\eqref{eq:240}
takes a solution to Eq.{~}\eqref{eq:020}
to the function verifying the same equation.
If  $\alpha =\half\pi $ then
this map of the space of solutions to Eq.{~}\eqref{eq:020}
into itself reduces to the identical map.

On the other hand, {\it in case of integer $\llll$}, there exist two
additional (as compared to the case of generic $\llll$)
transformations preserving the space of solutions 
to Eq.{~}\eqref{eq:220} \cite{BT1}.
One of them,
which we denote $\opB$,
can be represented as follows:
\begin{equation}
                    \label{eq:250}
\begin{aligned}
\opB\!\!: \somEE(z)\mapsto
\opB[\somEE\,](z)
\colonequals
&
(-1)^\llll \,  2\omega \,
z^{-\llll+1} e^{\mu(z+z^{-1})}
\left(z^2 \rrr(-z) \somEE'(-z)+ \sss(-z) \somEE(-z)\right).
\\
\end{aligned}
\end{equation}
The invariance
of the space of solutions to  Eq.{~}\eqref{eq:220}  with respect to $\opB$
can be established by straightforward
computations, provided the following property of
the polynomials $\ppp,\qqq,\rrr,\sss$
\begin{equation}
               \label{eq:260}
\begin{aligned}
\ppp(-z)
          =&\;
        (-1)^{\llll+1} (\lambda + \mu^2)^{-1} \left(\mu z^2 \rrr(z) + \sss(z)\right)
,
\\
\qqq(-z)
          =&\;
      \mu z^2 \ppp(z) + \qqq(z)
     +(-1)^\llll (\lambda + \mu^2)^{-1} \mu z^2 \left(\mu z^2 \rrr(z) + \sss(z)\right)
,
\\
\rrr(-z)
          =&\;
           \rrr(z)
,
\\
\sss(-z)
          =&\;
      (-1)^{\llll+1} (\lambda + \mu^2) \ppp(z)
    - \mu z^2 \rrr(z);
\end{aligned}
\end{equation}
and the differential equations
\begin{equation}
              \label{eq:270}
\begin{aligned}
z^2 \ppp'\argz=&\big(\mu + (\llll-1) z\big) \ppp\argz - \qqq\argz + (-1)^\llll z^2 \rrr\argz
,
 \\
\qqq'\argz=& \big(\lambda - (\llll+1) \mu z\big) \ppp\argz + \mu\, \qqq\argz + (-1)^\llll \sss\argz
,
 \\
z^2 \rrr'\argz=&(-1)^{\llll+1} \big(\lambda + \mu^2\big) \ppp\argz + z \big(2 (\llll-1) - \mu z\big) \rrr\argz - \sss\argz
,
 \\
z^2 \sss'\argz=&(-1)^{\llll+1} \big(\lambda + \mu^2\big) \qqq\argz
 + z^2 \big(\lambda - (\llll+1) \mu z\big) \rrr\argz + \big((\llll-1) z-\mu\big) \sss\argz
\end{aligned}
\end{equation}
which they obey \cite{BT1}
are taken into account.

Moreover, applying the operator $\opB$ twice and
utilizing the same reductions ensured by
\Eqs{~}\eqref{eq:260} and  \Eqs{~}\eqref{eq:270},
one finds that {\em on solutions to\/} Eq.{~}\eqref{eq:220},
the function-argument $\somEE$ is finally restored up to some constant factor and
{\em up to modification of its argument\/}
which undergoes, ultimately, the transformation
looking like a
full revolution around zero yielding no ultimate effect in
projection to
 $\mathbb{C}^*$
but
identical to the  \monodrom{}  transformation on the actual domain
\domain{}  of $\somEE$. This result
can be captured by means of  the following equality:
\begin{eqnarray}
           \label{eq:280}
\opB\circ\opB&=&
\delTa\cdot \mono.
\end{eqnarray}
In computation of the operator composition $\opB\circ\opB  $,
the factor $ \delTa $ appears 
originally
 in the following form
\begin{equation}
       \label{eq:290}
\delTa=z^{2(1-l)}\big(\ppp(z)\sss(z)-\qqq(z)\rrr(z)\big).
\end{equation}
However, a straightforward computation shows that $ \delTa $ is the first integral
of the system of differential equations \eqref{eq:270}
which the polynomials involved in its definition obey.
Thus $ \delTa $ does not actually depend on the
variable $z$ and can be determined %
 setting any value of the latter.
Substituting, in particular, $z=1$,
one obtains
\begin{equation}
           \label{eq:300}
\delTa=(2\omega)^{-2}\delTa_{+}\delTa_{-},
\end{equation}
where the factors on the right are defined in \eqref{eq:170}.

The equality 
\eqref{eq:280} and formulas \eqref{eq:170}
now say us the following. 
\begin{proposition}
  If
\begin{equation}
            \label{eq:310}
\delTa_{+}\not=0\not=\delTa_{-} \mbox{ or, equivalently, }
\ppp(1)^2\not=(2\omega)^2\rrr(1)^2
\end{equation}
then the linear operator $\opB$ \eqref{eq:250}
determines the automorphism of the space of solutions to Eq.{~}\eqref{eq:220}.
\end{proposition}
The violation of the condition \eqref{eq:310}
would impose  severe 
restrictions on the constant parameters involved in
Eq.{~}\eqref{eq:220}.
We assume to consider a generic case 
claiming
 \eqref{eq:310} to be fulfilled throughout.

The linear operator $\opB  $ acting on the two-dimensional linear space of solutions
to Eq.{~}\eqref{eq:220} can be presented with respect to any basis of this space as
some $2\times2$ matrix.
In particular,
it can be shown that in the basis $(\EEpm{+},\EEpm{-})$ introduced
above the matrix form of the operator $\opB  $ reads
\begin{equation}
           \label{eq:320}
\begin{aligned}
\mathbf{B}=
&
\Imi^\llll(2\omega)^{-1}e^{\mhalf P(0)}\times
\\&
\diag\left(
-e^{\mquarti\pi}\delTa_{+}\cos\big(\half(\varphi(0)-\half\pi )\big)^{-1}
,e^{\quarti\pi}\delTa_{-}\cos\big(\half(\varphi(0)+\half\pi )\big)^{-1}
\right)\times
\\&\hphantom{\diag}
\begin{pmatrix}
e^{\half{}P(\half T)}\CCC-e^{\half{}P(\mhalf T)}\DDD
&
\Imi(
e^{\half{}P(\half T)}\CCC+e^{\half{}P(\mhalf T)}\DDD
)
\\
\Imi(
e^{\half{}P(\half T)}\AAA+e^{\half{}P(\mhalf T)}\BBB
)
&
-e^{\half{}P(\half T)}\AAA+e^{\half{}P(\mhalf T)}\BBB
\end{pmatrix}.
\end{aligned}
\end{equation}

Let us note now that the numerator of the fraction
in Eq.{~}\eqref{eq:240} is a solution to
Eq.{~}\eqref{eq:220}
and the denominator is also a solution  get
with the modified argument ($1/z$ substituted in place of $z$),
the former and the later being
mutually
conjugated since
$\overline{\EEpm{\pm}{(\bar z)}}=\EEpm{\pm}(z)$
under the conditions assumed.

If one replaces in \eqref{eq:240}, formally, the functions $\EEpm{\pm}$
in the numerator
with the functions $\opB\EEpm{\pm}$ expanding  them further as  linear combinations
of the
original  $\EEpm{\pm}$ derived
making use of the matrix \eqref{eq:320},
and carry out the corresponding transformation
of the denominator preserving its complex conjugacy with the numerator
on the unit circle,
then the formula similar to \eqref{eq:240}
but with distinct parameter $\alpha$ results.
Repeating such a transformation twice, one comes, in view of \eqref{eq:280},
to the original
functions $\EEpm{\pm}$ with the original coefficients $ \cos(\half\alpha) $,
 $ \sin(\half\alpha) $ (times the constant factor of
$\delTa$ which cancels out) but with
arguments undergone the \monodrom{}  transformation.
In other words, the transformation of
\eqref{eq:240}
induced by the map \eqref{eq:250},
repeated twice, results in the \monodrom{}  transformation of the function
$ {\eiphi}^{(\alpha)}$.
Setting  $\alpha=\half\pi$,
performing such a transformation of $ {\eiphi}^{(\alpha)}$
once, and eliminating
the functions $\EEpm{\pm}$ by means of their expansions \eqref{eq:200},
the formula \eqref{eq:130} results.

The functions $\calF$ and $\CCcalF$
(see \Eqs{~}\eqref{eq:150}, \eqref{eq:160})
which
have been used for determination of the function $\eP$
alone (see Eq.{~}\eqref{eq:140})
are actually of notable   interest  in their own rights.
Such functions are closely related to solutions to \eqref{eq:020}.
They can be defined as
solutions to the linear differential equations
\begin{equation}
            \label{eq:330}
2\,\Imi\,\omega\,z\,\calF'=- {\eiphi}(\calF- \CCcalF ),\;
2\,\Imi\,\omega\,z\,\CCcalF'= {\eiphi}^{-1}(\calF- \CCcalF )
\end{equation}
obeying the initial conditions
\begin{equation}
            \label{eq:340}
\calF(1)=\Imi, \;\CCcalF(1)=-\Imi.
\end{equation}
If these are fulfilled,
then the difference $\eP=(2\Imi)^{-1} (\calF- \CCcalF ) $ (see  Eq.{~}\eqref{eq:140})
obeys Eq.{~}\eqref{eq:080}
and represents
 the analytic continuation of the function
$e^{P(t)}=\exp\int_0^t\cos\varphi(\tilde t)\,d \tilde t  $,
where $2\cos\varphi(t)={\eiphi}(e^{\Imi\omega t})+{\eiphi}(e^{\Imi\omega t})^{-1}$,
from an arc of the unit circle to its vicinity in $\mathbb{C}^*$.

Concerning the very formulas \eqref{eq:150}, \eqref{eq:160},
the fulfillment of \Eqs{~}\eqref{eq:330}, \eqref{eq:340}
by the functions the former define 
for the corresponding right hand side factors ${\eiphi}_B^{\pm1}$
defined by  Eq.{~}\eqref{eq:130}
is verified by straightforward computations.

In conclusion, it should be emphasized
that the transformation \eqref{eq:150},
taking solutions to the nonlinear equation \eqref{eq:020}
to solutions of the same equation
(and, then, determining the associated transformation on the
space of solutions to Eq.{~}\eqref{eq:010}),
arises here 
as a byproduct
of the
specific
symmetry of the space of solutions to the linear equation
\eqref{eq:220}.
Such a symmetry
has been shown to exists
in case of integer value of the
parameter $\llll$ (sometimes called the order).
The existence of analogue of the above $B$-transformation
under less restrictive conditions  remains an open problem.

\end{document}